\documentclass{article}

\usepackage{mwe} 
\usepackage{caption}
\usepackage{amsfonts,mathrsfs}
\usepackage{amsmath, amsthm, amssymb}

\usepackage{float}
\usepackage{subfigure}
\newfloat{figure}{H}{lof}
\newfloat{table}{H}{lot}
\floatname{figure}{\figurename}
\floatname{table}{\tablename}

\usepackage{mathtools}
\usepackage{etoolbox}
\usepackage[colorlinks=true]{hyperref}
\usepackage{hyperref}
\hypersetup{linkcolor=blue}

\usepackage{enumerate}

\usepackage{tcolorbox}
\usepackage{colortbl}
\usepackage{bbm}
\usepackage{nicematrix}

\usepackage[latin1]{inputenc}

\usepackage{float}
\newfloat{figure}{H}{lof}
\newfloat{table}{H}{lot}
\floatname{figure}{\figurename}
\floatname{table}{\tablename}

\numberwithin{equation}{section}

\numberwithin{equation}{section}

\newcommand {\CR}{{\rm CR}}

\begin{document}

\title{\textbf{Forecasting  the changes between endemic and epidemic phases of a contagious disease, with the example of COVID-19}}
\author{\textsc{Jacques Demongeot $^{1}$, Pierre Magal$^{2}$\footnote{Corresponding author. e-mail: \href{mailto:pierre.magal@u-bordeaux.fr}{pierre.magal@u-bordeaux.fr}}, and Kayode Oshnubi $^{1}$}\\
	{\small \textit{$^{1}$ AGEIS laboratory, UGA, 38700 La Tronche, France}}\\
	{\small \textit{$^{2}$ Univ. Bordeaux, IMB, UMR 5251, F-33400 Talence, France.}} \\
	{\small \textit{CNRS, IMB, UMR 5251, F-33400 Talence, France.}}\\
}
\maketitle

\begin{abstract}
Predicting the endemic/epidemic transition during the temporal evolution of a contagious disease. \\
\textit{Methods: } Defining indicators for detecting the transition endemic/epidemic, with four scalars to be compared, calculated from the daily reported news cases: coefficient of variation, skewness, kurtosis, and entropy. The indicators selected are related to the shape of the empirical distribution of the new cases observed over 14 days. This duration has been chosen to smooth out the effect of weekends when fewer new cases are registered. For finding a forecasting variable, we have used the PCA (principal component analysis), whose first principal component (a linear combination of the selected indicators) explains a large part of the observed variance and can then be used as a predictor of the phenomenon studied (here the occurrence of an epidemic wave).\\
\textit{Results:} A score has been built from the four proposed indicators using a Principal Component Analysis (PCA), which allows an acceptable level of forecasting performance by giving a realistic retro-predicted date for the rupture of the stationary endemic model corresponding to the entrance in the epidemic exponential growth phase. This score is applied to the retro-prediction of the limits of the different phases of the COVID-19 outbreak in successive endemic/epidemic transitions and three countries, France, India, and  Japan. \\
\textit{Conclusion: }We provided a new forecasting method for predicting an epidemic wave occurring after an endemic phase for a contagious disease.
\end{abstract}

\medskip


\vspace{0.5cm}
\noindent  \textit{The paper is dedicated to James D. Murray, whose pioneering work in mathematical biology we admire.}

\section{Introduction}

Finding a reliable prediction method of the frontiers between different stationary and non-stationary periods of a time series is a challenging problem. Since the seminal work by Deshayes and Picard on the stationarity rupture in time series \cite{Deshayes-Picard1979, Deshayes-Picard1984, Picard1985}, many works have dealt with the break in stationarity  \cite{4, 5, 6, 7, 8, 9}, the most recent using the concepts of functional statistics \cite{10, 11, 12, 13, 14, 15, 16}. Indeed, stationarity is crucial as many forecasting models of time series rely on stationarity for easy modeling and obtaining reliable results. A stationary time series presents statistical properties which do not change over time, as the empirical distribution of the random variable observed in the series, with its main characteristic parameters, mean, coefficient of variation, moments, and entropy.
In the event of a break in stationarity, there may be a sudden transition with a sudden change in the values of these parameters and the appearance of a non-constant trend. The problem of the existence of this transition arises with particular acuity in the case of contagious diseases, which alternate stationary endemic periods and epidemic peaks having an initial exponential trend, which must be predicted to prevent the spread of the disease does not give rise to a pandemic.

The term endemic phase is understood to mean a period in which there is an equilibrium in a model whose parameters have changed in value following an epidemic phase, due to mitigation measures or a change in the virulence of the infectious agent. In the case of chronic diseases observed outside epidemic phases (bacterial meningitis, rabies, smallpox before its eradication, etc.), the definition of endemicity corresponds more to the sporadic appearance (by birth or human displacement) of individuals much more susceptible than the general population to the infectious agent \cite{Bernoulli1760,Blower,Murray}.

We will propose in this article a method to estimate the breakdown of endemic stationarity based on four parameters linked to the empirical distribution of the number of daily reported new  cases of COVID-19 in several countries, parameters whose isolated or joint predictive power will be analyzed. These parameters are the coefficient of variation , the skewness, the kurtosis,  and the entropy of the stationary empirical measure calculated in a moving window.

The pair represented by the succession of an endemic phase and an epidemic wave in COVID-19 outbreak  can be considered as a functional unit, the break between the two phases having to be found \cite{16_1, 16_2, 16_3, 16_4, 16_5, 16_6}. The endemic phase is characterized by a low average level of new cases, with low variance. At the start of the epidemic phase, the average number of cases will grow exponentially, and the standard deviation will grow proportionally at the beginning, then saturate, like that of an additive noise in part independent of growth, which explains the increase then decrease in the coefficient of variation and kurtosis, therefore those of the first principal component at the endemic/epidemic boundary.

\section{Materials and Methods}

We use in the following a moving window of length 14 days for calculating the empirical distribution of the random variable equal to the number of daily reported new cases. The empirical distribution $N_t$ on day $t$ is obtained from the daily number of reported new cases considered as a random variable $N_t= \left(N(t-13), N(t-12), \ldots, N(t)\right)$.

The length of the window has been chosen for eliminating the effect of the week periodicity observed in data due to the lack of reporting each weekend. Indeed, daily new infection cases are highly affected by weekends, such that new case numbers are lowest at the start of the week, and increase afterwards \cite{16A}. 
\subsection{Empirical first four moments}
In the following we will use the terminology \textit{endemic period}, to describe a period during which the daily new cases occurs randomly around some mean value. An \textit{epidemic wave} is a period during which the  daily new cases occurs by contact between susceptible and infected. 

Our goal in this paper is to explore the transition between the endemic period and  the following epidemic wave which will be studied by calculating several parameters in a moving window around the frontier of 14 days on which we suspect this transition occurred. 

\medskip 
We consider the first four moments of $N_t$. We start with the \textit{mean}  
\begin{equation}
	\mu=E(N_t) = \dfrac{ \displaystyle \sum_{i=0}^{13} N(t-i)}{14},
\end{equation}
where $E$ is the expectation operator, with the \textit{standard deviation}
\begin{equation}
	\sigma=E\left(\left(N_t-\mu\right)^2\right)^{1/2} = \sqrt{ \dfrac{ \displaystyle \sum_{i=0}^{13} \left(N(t-i)-\mu\right)^2 }{14}}.
\end{equation}

From these two first parameters, we can compute the \textit{coefficient of variation} 
\begin{equation}
CV(N_t)=\dfrac{\sigma}{\mu}. 
\end{equation}                                                                                     
The \textit{skewness} of the random variable $N_t$  is the third standardized moment, defined as
\begin{equation}
		Skew(N_t)=E\left(\left(\dfrac{N_t-\mu}{\sigma} \right)^3\right).
\end{equation}
Recall that the skewness verifies 
\begin{equation}
	Skew(N_t)=   \dfrac{E(N_t^3)-3 \mu \sigma^2 - \mu^3}{\sigma^3}= \dfrac{E(N_t^3)}{\sigma^3}- \left(3\dfrac{1}{CV}+\dfrac{1}{CV^3} \right). 
\end{equation} 
The \textit{kurtosis} is the fourth standardized moment, defined as 
\begin{equation}
	Kurt(N_t)=E\left(\left(\dfrac{N_t-\mu}{\sigma} \right)^4\right) .
\end{equation}
The \textit{empirical entropy} $\mathcal{E}$ of the empirical distribution is defined as follows:
\begin{equation}
\mathcal{E} (N_t) = - \sum_{i=1:d \text{ with } p_i >0} p_i   \log p_i,         
\end{equation}
where the $p_i$ are the weights of a histogram on d value intervals of $N_t$. In the Results' section, we use the \textit{approximate entropy}. We refer to \cite{Pincus} for more details.

\subsection{Phenomenological model used for multiple epidemic waves}
To represent the data, we used a phenomenological model to fit the curve of cumulative reported cases. Such an idea is not new since it was already proposed by Bernoulli \cite{Bernoulli1760} in 1760 in the context of the smallpox epidemic. Here we used the so-called Bernoulli--Verhulst \cite{Verhulst1838} model to describe the epidemic phase.  Bernoulli \cite{Bernoulli1760} investigated an epidemic phase followed by an endemic phase. This appears clearly in Figures 9 and 10 of the paper by Dietz, and Heesterbeek \cite{Dietz-Heesterbeek}  who revisited the original article of Bernoulli. We also refer to Blower \cite{Blower} for another article revisiting the original work of Bernoulli. Several works comparing cumulative reported cases data and the Bernoulli--Verhulst model appear in the literature  (see \cite{Hsieh, Wang-Wu-Yang, Zhou-Yan}).  The Bernoulli--Verhulst model is sometimes called Richard's model,\index{Richardsmodel@Richard's model} although Richard's work came much later in 1959.

The phenomenological model deals with data series of new infectious cases decomposed into two successive phases: 1) endemic phases followed by 2) epidemic phases.

\medskip
\noindent \textbf{Endemic phase:}\index{endemic phase} During the endemic phase, the dynamics of new cases appears to fluctuate around an average value independently of the number of cases. Therefore the average cumulative number of  cases is given by
\begin{equation} \label{15.11}
	\CR(t)=N_0+ (t-t_{0}) \times a, \text{ for } t \in [t_0 , t_1],
\end{equation}
where $t_0$ {denotes} the beginning of the endemic phase, $N_0$ is the number of new cases at time $t_0$, and $a$ is the average value of the daily number of new cases.

We assume that the average daily number of {new} cases is constant. Therefore the daily number of {new} cases is given by
\begin{equation} \label{15.12}
	\CR'(t)=  a.
\end{equation}
\medskip
\noindent \textbf{Epidemic phase:}\index{endemic phase} In the epidemic phase, the new cases  are contributing to produce  secondary cases. Therefore the daily number of new cases is no longer constant, but varies with time as follows
\begin{equation}\label{15.13}
	\CR(t)=	N_{\rm{base}}+ \dfrac{\mathrm{e}^{\chi (t-t_0)} N_0 }{\left[  1+ \dfrac{ N_0^\theta}{N_\infty^\theta}  \left(\mathrm{e}^{\chi \theta   \left(t- t_0 \right)  } -1 \right) \right]^{ 1/\theta} }, \text{ for } t \in [t_0 , t_1].
\end{equation}
In other words, the daily number of {new} cases  follows the Bernoulli--Verhulst \cite{Bernoulli1760, Verhulst1838} equation. Namely, by setting
\begin{equation}  \label{15.14}
	N(t)=\CR(t)-N_{\rm{base}},
\end{equation} we obtain
\begin{equation} \label{15.15}
	N'(t)= \chi \,  N(t) \, \left[ 1 - \left( \dfrac{N(t) }{N_\infty}\right)^\theta \right],
\end{equation}
completed with the initial value
\begin{equation*}
	N(t_0)=N_0.
\end{equation*}
In the model, $N_{\rm{base}}+N_0$  corresponds to the value $\CR(t_0)$ of the cumulative number of cases at time $t=t_0$. The parameter $N_\infty+N_{\rm{base}}$ is the maximal value of the cumulative reported {cases} after the time $t=t_0$.   $\chi>0$ is a Malthusian growth  parameter, and $\theta$ {regulates} the speed at which $\CR(t)$ {increases} to $N_\infty+N_{\rm{base}}$.

\section{Results}

Here we use  cumulative numbers of reported new cases for  COVID-19 in France, India, and Japan  taken from WHO \cite{Data}. Data shows a succession of endemic periods (yellow background color regions) followed by epidemic waves (blue background color regions).
\subsection{Data for France}

\begin{figure}[H]
		\begin{center}
				\includegraphics[scale=0.25]{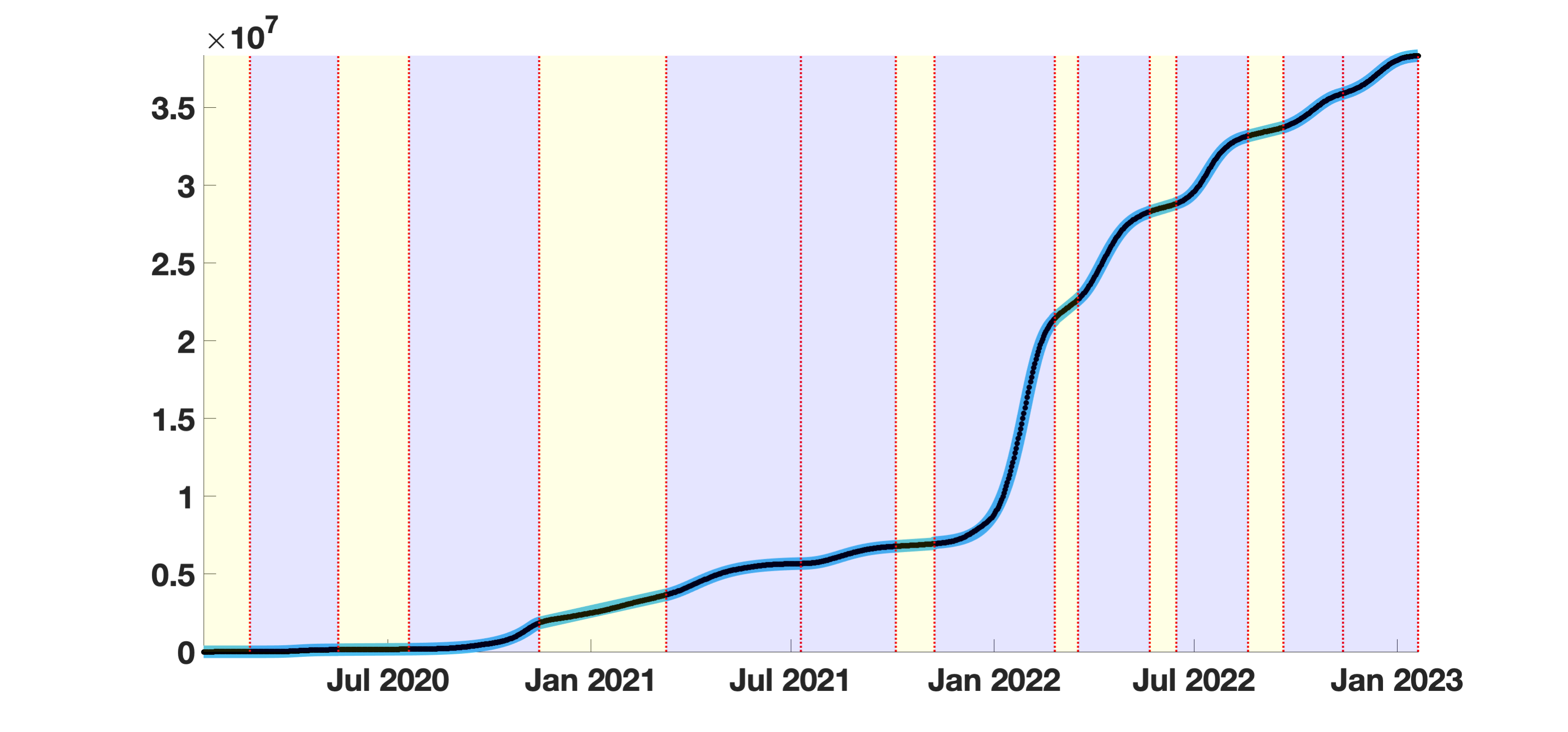}
			\end{center}
		\caption{\textit{In this figure we plot   in blue the phenomenological model  and in black the data.  Data is the cumulative  reported  number  of  new cases  with a 14-day rolling average.  }}\label{Fig1}
	\end{figure}
\begin{figure}[H]
	\begin{center}
		\includegraphics[scale=0.25]{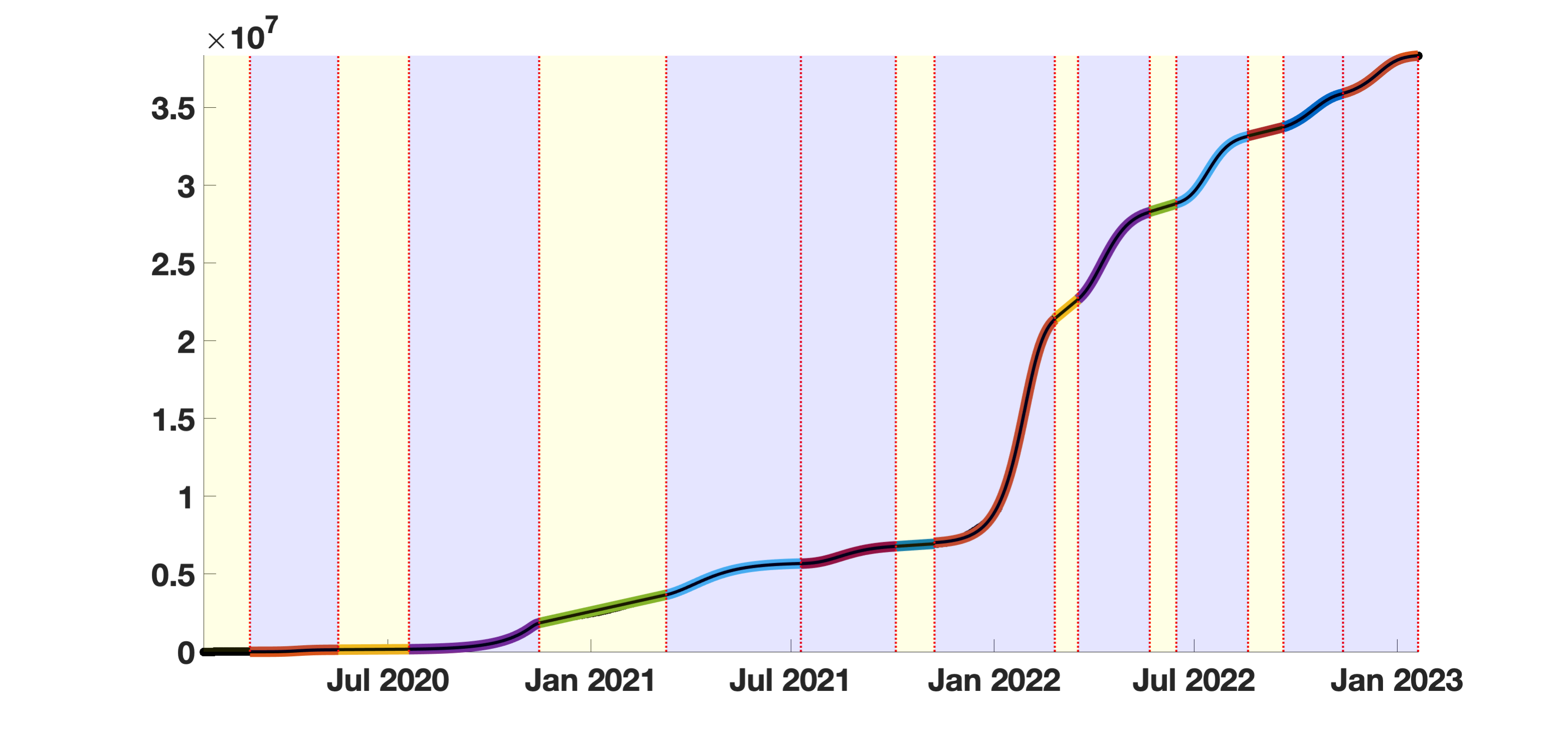}
	\end{center}
	\caption{\textit{In this figure we plot  with multiple colors the phenomenological models  obtained for each period.   }}\label{Fig2}
\end{figure}

In Figure \ref{Fig2}, each colored segment corresponds to a single endemic or epidemic phenomenological model. This change of color may occur due to a change of dynamic inside an epidemic phase when a second wave comes before the end of the previous one.

\begin{figure}[H]
	\begin{center}
		\includegraphics[scale=0.25]{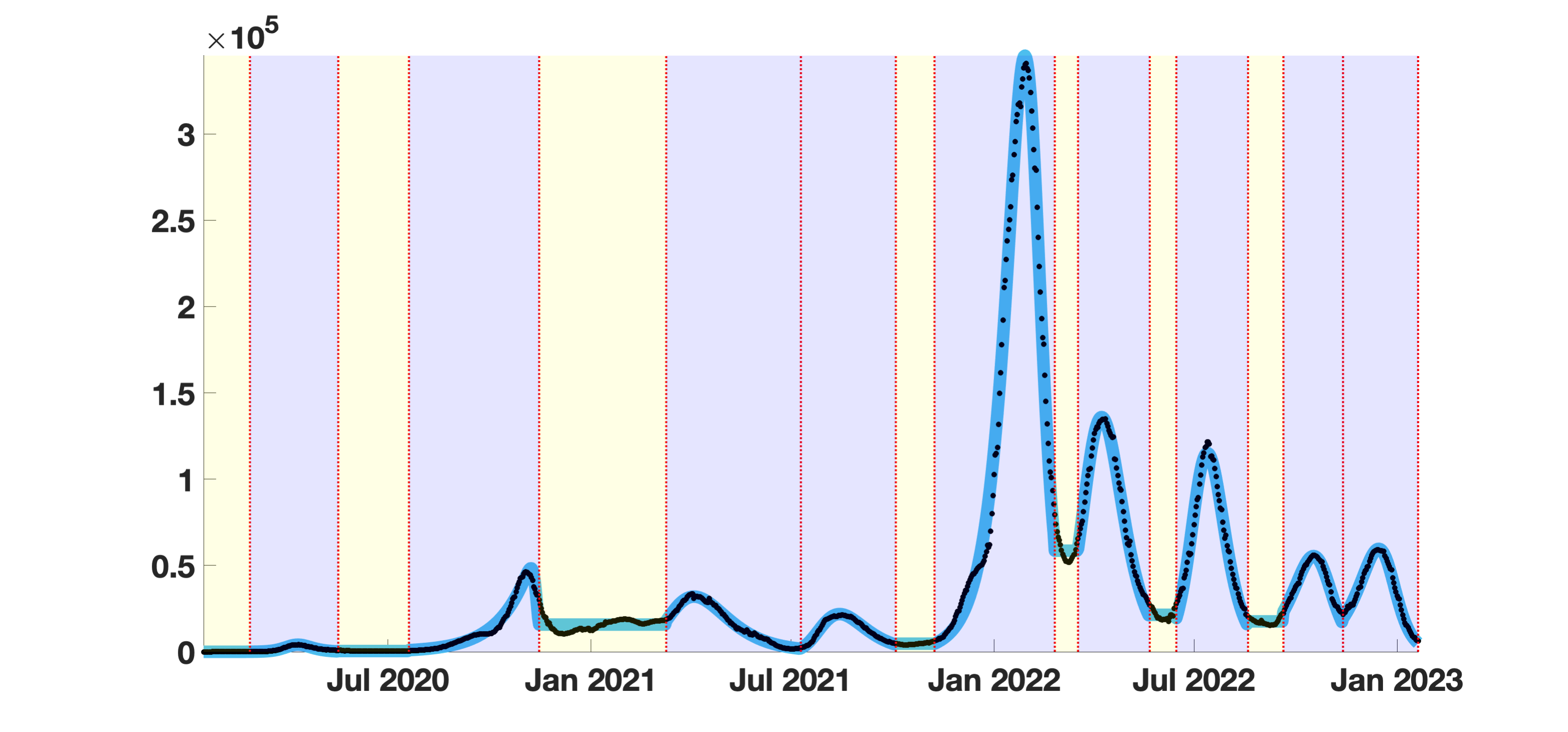}
	\end{center}
	\caption{\textit{In this figure we plot   in blue the first derivative of the phenomenological model  and in black the data. Data is the daily  reported  number of  new cases  with a 14-day rolling average.    }}\label{Fig3}
\end{figure}

By making a principal component analysis (i.e. the matlab function PCA) between the standardized variables $CV_s(N_t), Skew_s(N_t),  Kurt_s(N_t),
 \mathcal{E}_s(N_t)$, we obtain the percentage of the variance explained  by each principal component
$$
Explain= \left( \begin{array}{cccc}
70.71\\
21.92\\
4.94\\
2.43\\
\end{array}
\right)
$$ 
and the matrix giving the projection coefficients of the principal components 
$$
coeff =\left( \begin{array}{cccc}

 0.5527 &  -0.1480 &   0.7163 &  -0.3995\\
0.5631  & -0.2162  & -0.0348    & 0.7968\\
0.5577   & -0.0795 &  -0.6955 &  -0.4461\\
0.2577  &  0.9618   & 0.0449  &  0.0808
\end{array}
\right).
$$  
By using the first column of the above matrix, we deduce the first principal component 
\begin{equation} \label{8}
0.55 CV_s(N_t)+  0.56  Skew_s(N_t)+ 0.56 Kurt_s(N_t)+ 0.26 \mathcal{E}_s(N_t)	
\end{equation}
which explains $70.71 \%$ of the variability.

We deduce that Kurtosis, Skewness, and the coefficient of variation (in decreasing order of importance) best explain the variability.

\begin{figure}[H]
	\begin{center}
		\includegraphics[scale=0.25]{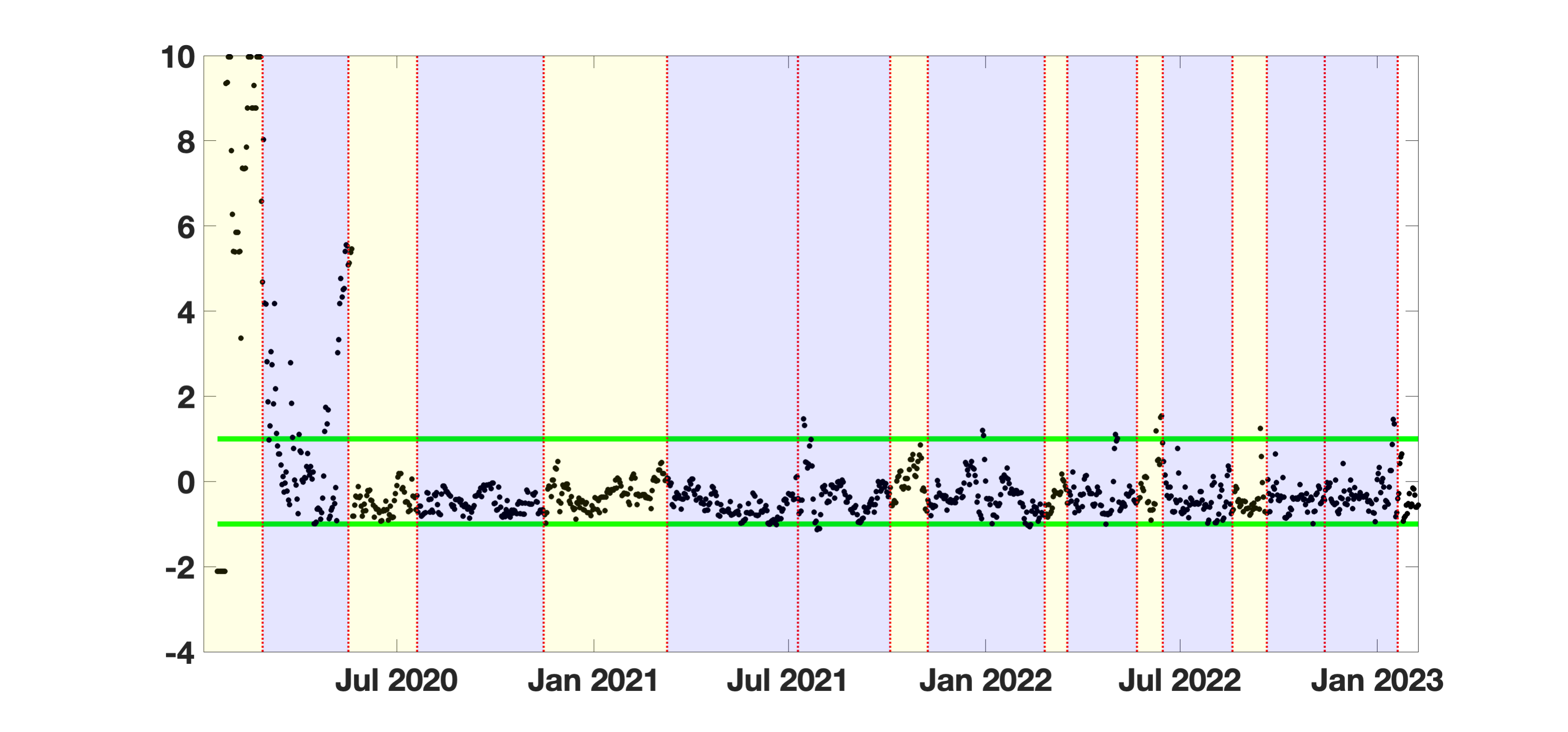}
	\end{center}
	\caption{\textit{In this figure we plot  the first principal component for France (see formula \eqref{8}).   The horizontal green lines correspond to $\pm 1$.  }}\label{Fig4}
\end{figure}

\subsection{Data for India}
In this subsection, we consider the data for India. We present the same curves as for France, describing successively raw data and simulated results by the phenomenological models.

\begin{figure}[H]
	\begin{center}
		\includegraphics[scale=0.25]{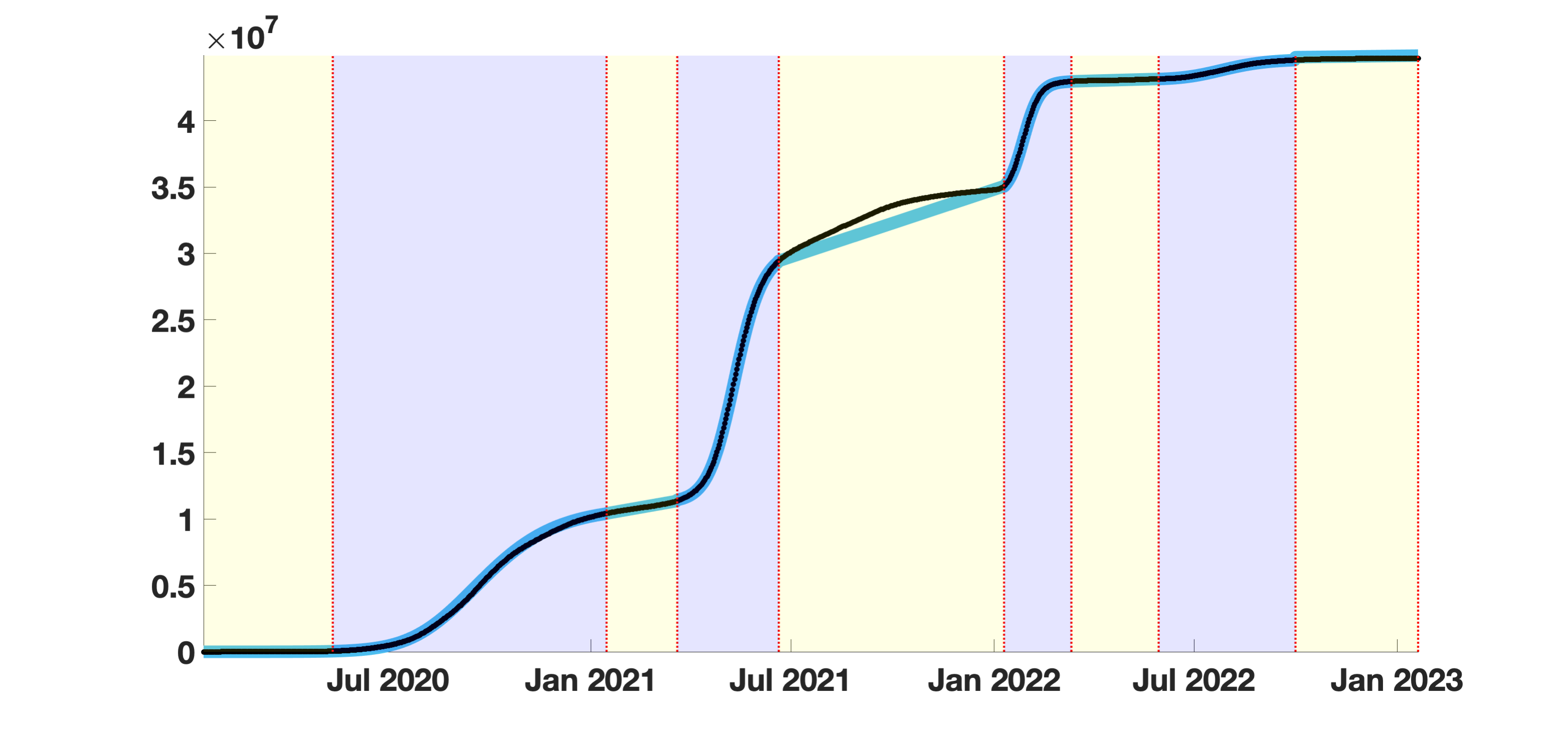}
	\end{center}
	\caption{\textit{In this figure we plot   in blue the phenomenological model  and in black the data.  Data is the cumulative  reported  number  of  new cases  with a 14-day rolling average.  }}\label{Fig5}
\end{figure}
\begin{figure}[H]
	\begin{center}
		\includegraphics[scale=0.25]{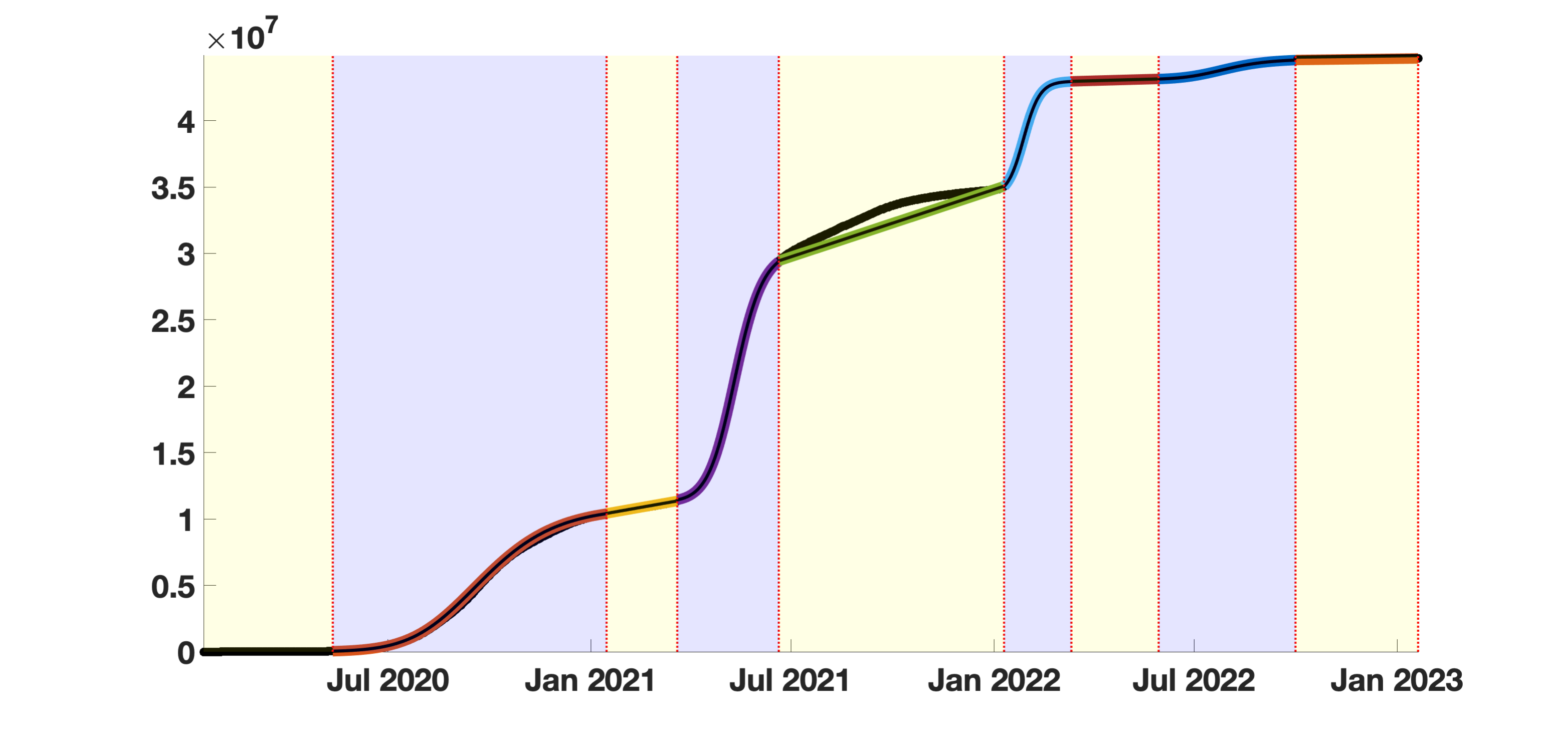}
	\end{center}
	\caption{\textit{In this figure we plot  with multiple colors the phenomenological models  obtained for each period.   }}\label{Fig6}
\end{figure}

In Figure \ref{Fig6}, each colored segment corresponds to a single endemic or epidemic phenomenological model. This change of color may occur due to a change of dynamic inside an epidemic phase when a second wave comes before the end of the previous one. 
\begin{figure}[H]
	\begin{center}
		\includegraphics[scale=0.25]{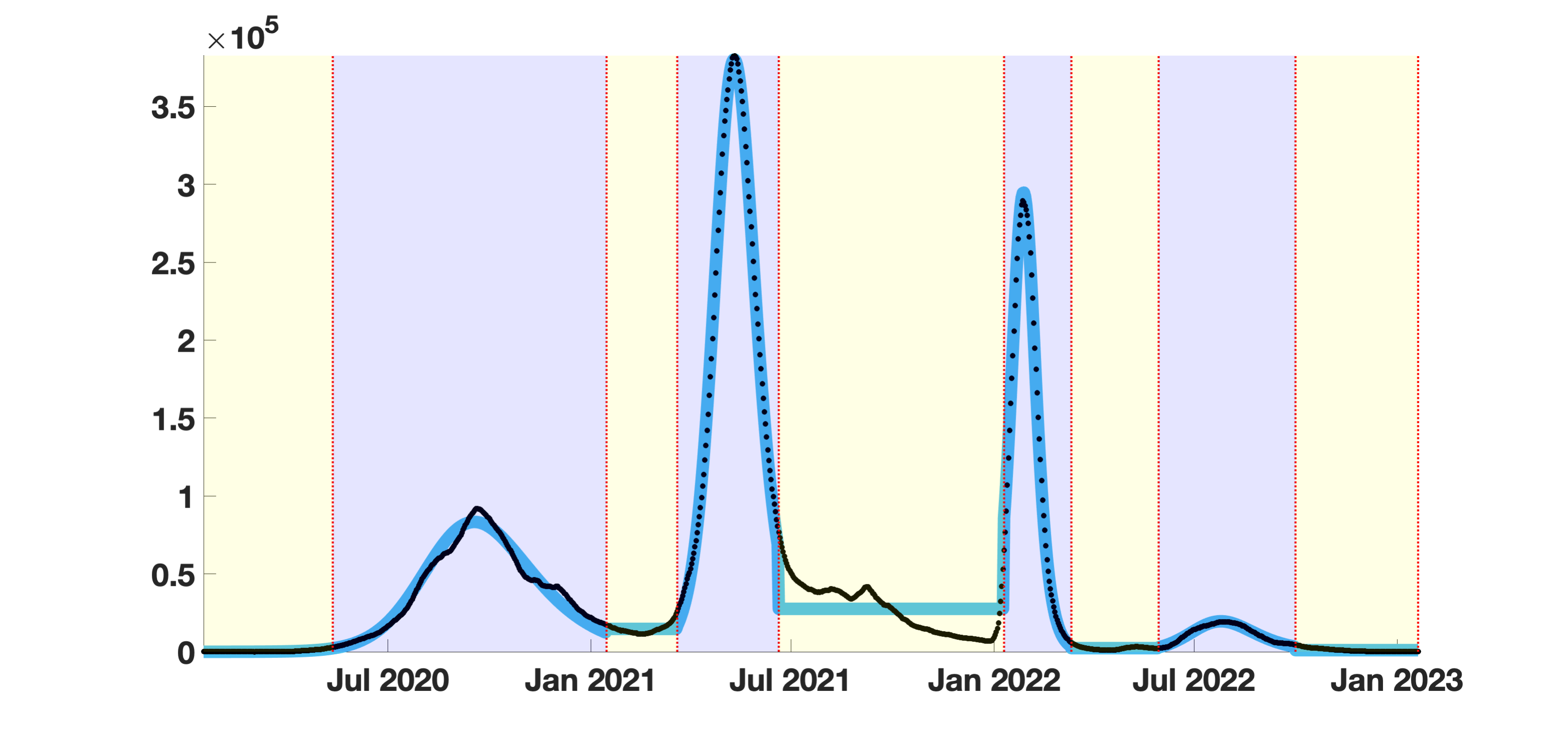}
	\end{center}
	\caption{\textit{In this figure we plot   in blue the first derivative of the phenomenological model  and in black the data. Data is the daily  reported  number of  new cases  with a 14-day rolling average.    }}\label{Fig7}
\end{figure}
By making a principal component analysis (i.e. the matlab function PCA) between the standardized variables $CV_s(N_t), Skew_s(N_t),  Kurt_s(N_t)$, and $ \mathcal{E}_s(N_t)$, we obtain the percentage of the variance explained  by each principal component
$$
Explain= \left( \begin{array}{cccc}
53.13\\
22.55\\
14.22\\
10.09 \\
\end{array}
\right)
$$ 
and the matrix giving the projection coefficients of the principal components 
$$
coeff =\left( \begin{array}{cccc}
  0.5161 &  -0.2296 &   0.8212  & -0.0810\\
0.5714  & -0.0587 &  -0.4434 &  -0.6881\\
0.5660   &-0.2235  &-0.3478 &   0.7133\\
0.2946  &  0.9455  &  0.0896 &   0.1062
\end{array}
\right).
$$  
By using the first column of the above matrix, we deduce the first principal component 
\begin{equation} \label{9}
	0.52 CV_s(N_t)+  0.57  Skew_s(N_t)+ 0.57 Kurt_s(N_t)+ 0.29 \mathcal{E}_s(N_t)	
\end{equation}
which explains $53 \%$ of the variability.

\begin{figure}[H]
	\begin{center}
		
		\includegraphics[scale=0.25]{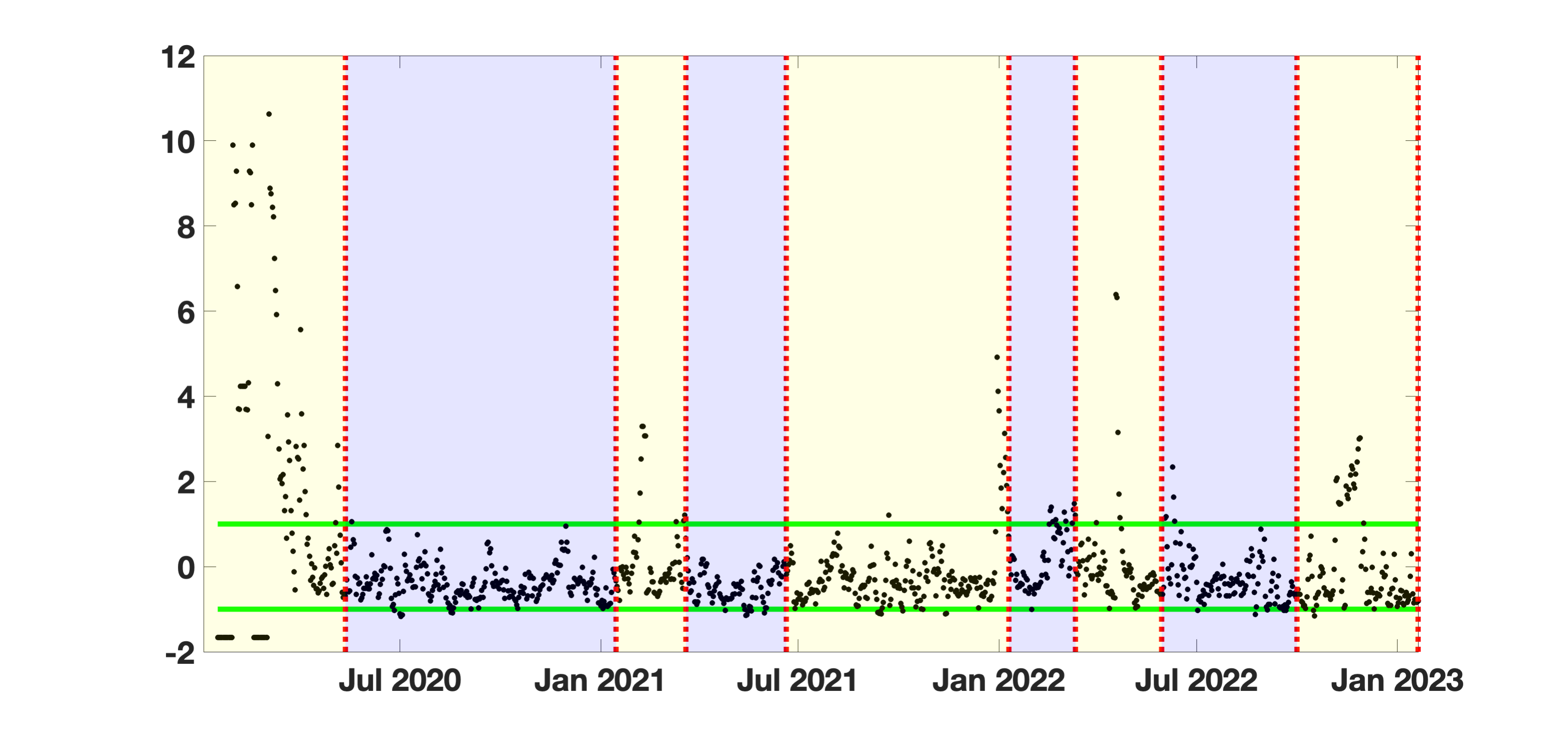}
		
	\end{center}
	\caption{\textit{In this figure we plot  the first principal component for India (see formula \eqref{9}).  The horizontal green lines correspond to the values $\pm 1$.   }}\label{Fig8}
\end{figure}

\subsection{Data for Japan}
In this subsection, we consider the data for Japan. We present the same curves as for France and India, describing successively raw data and simulated results by the phenomenological models.

\begin{figure}[H]
	\begin{center}
		\includegraphics[scale=0.25]{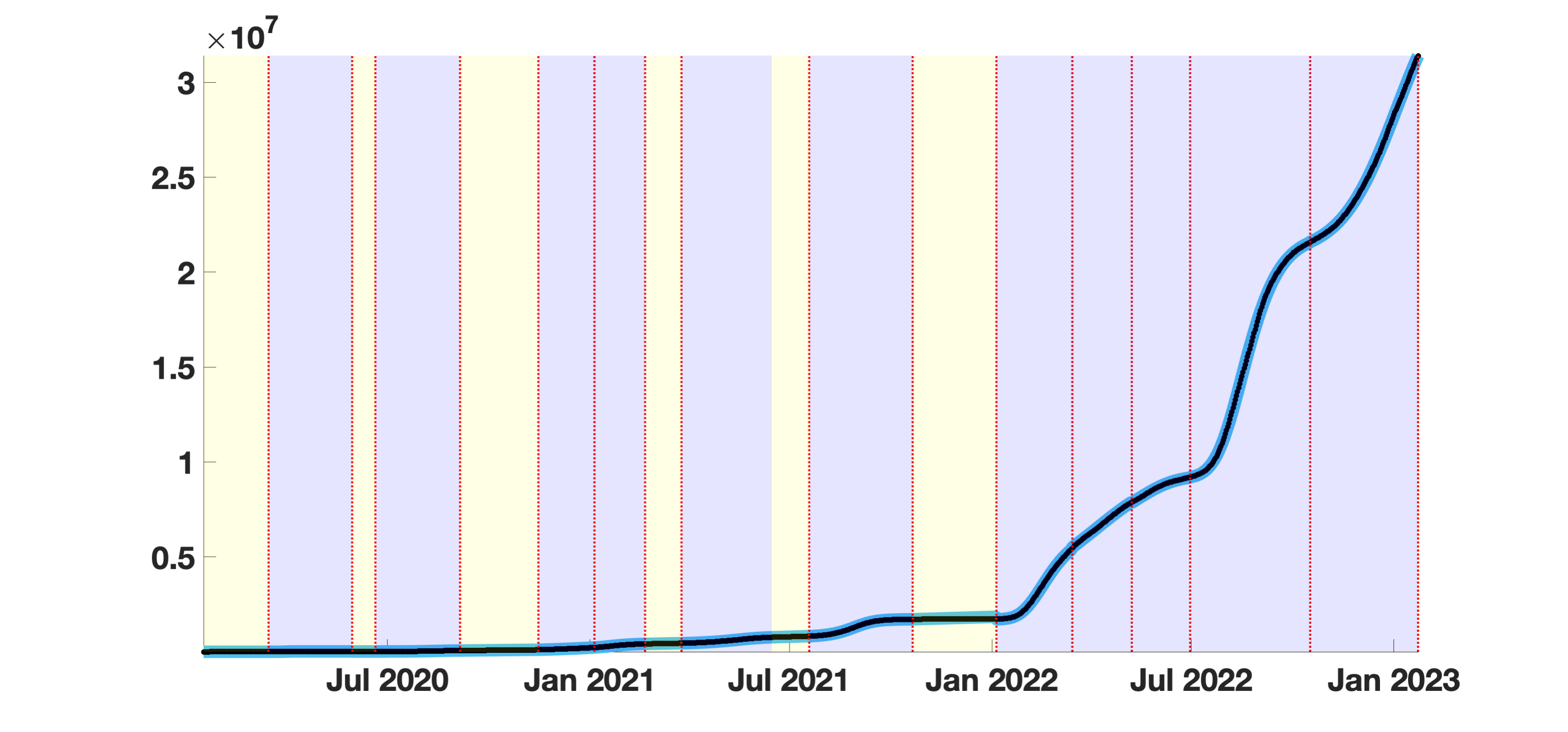}
	\end{center}
	\caption{\textit{In this figure we plot   in blue the phenomenological model  and in black the data.  Data is the cumulative  reported  number  of  new cases  with a 14-day rolling average.   }}\label{Fig9}
\end{figure}

\begin{figure}[H]
	\begin{center}
		\includegraphics[scale=0.25]{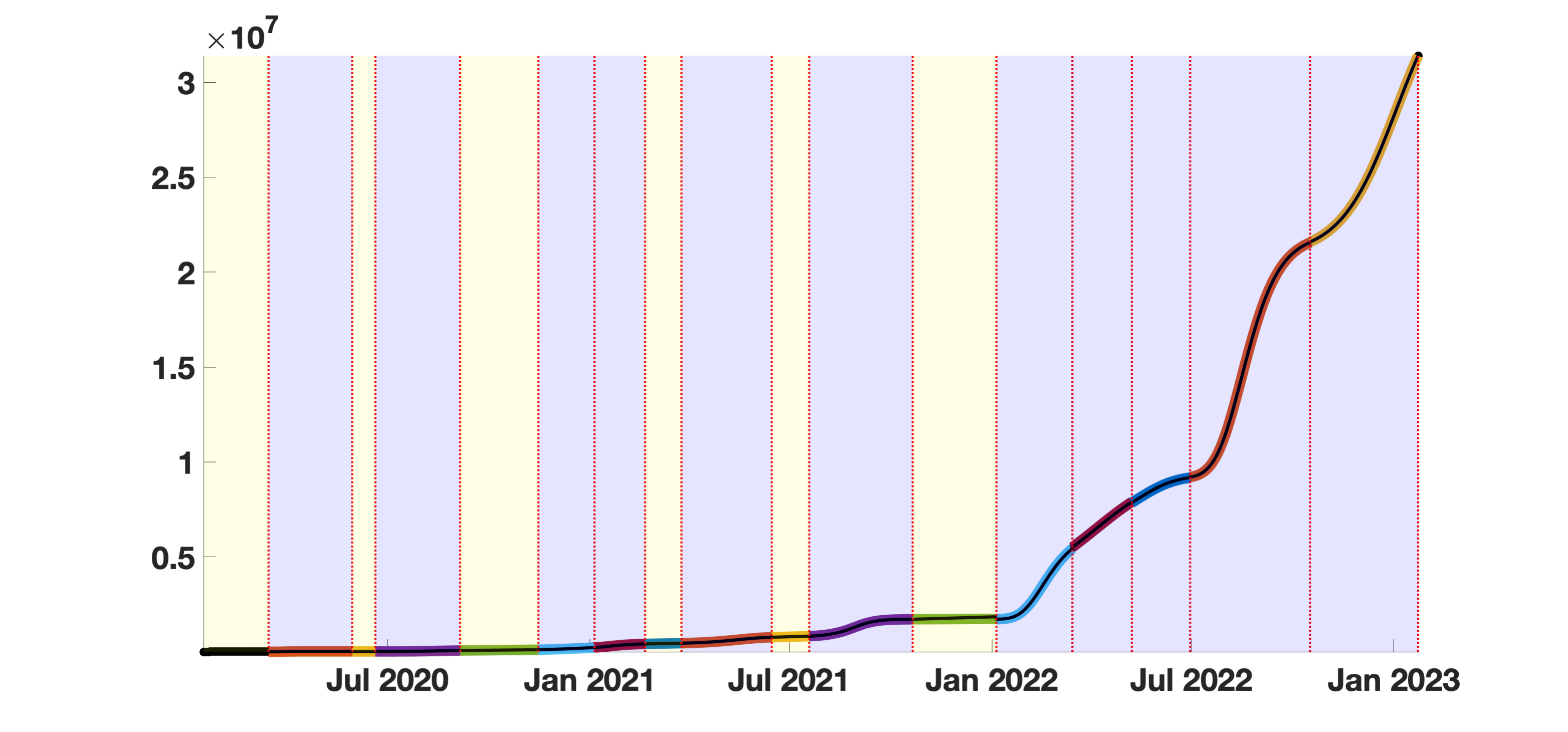}
	\end{center}
	\caption{\textit{In this figure we plot  with multiple colors the phenomenological models  obtained for each period.    }}\label{Fig10}
\end{figure}


\begin{figure}[H]
	\begin{center}
		\includegraphics[scale=0.25]{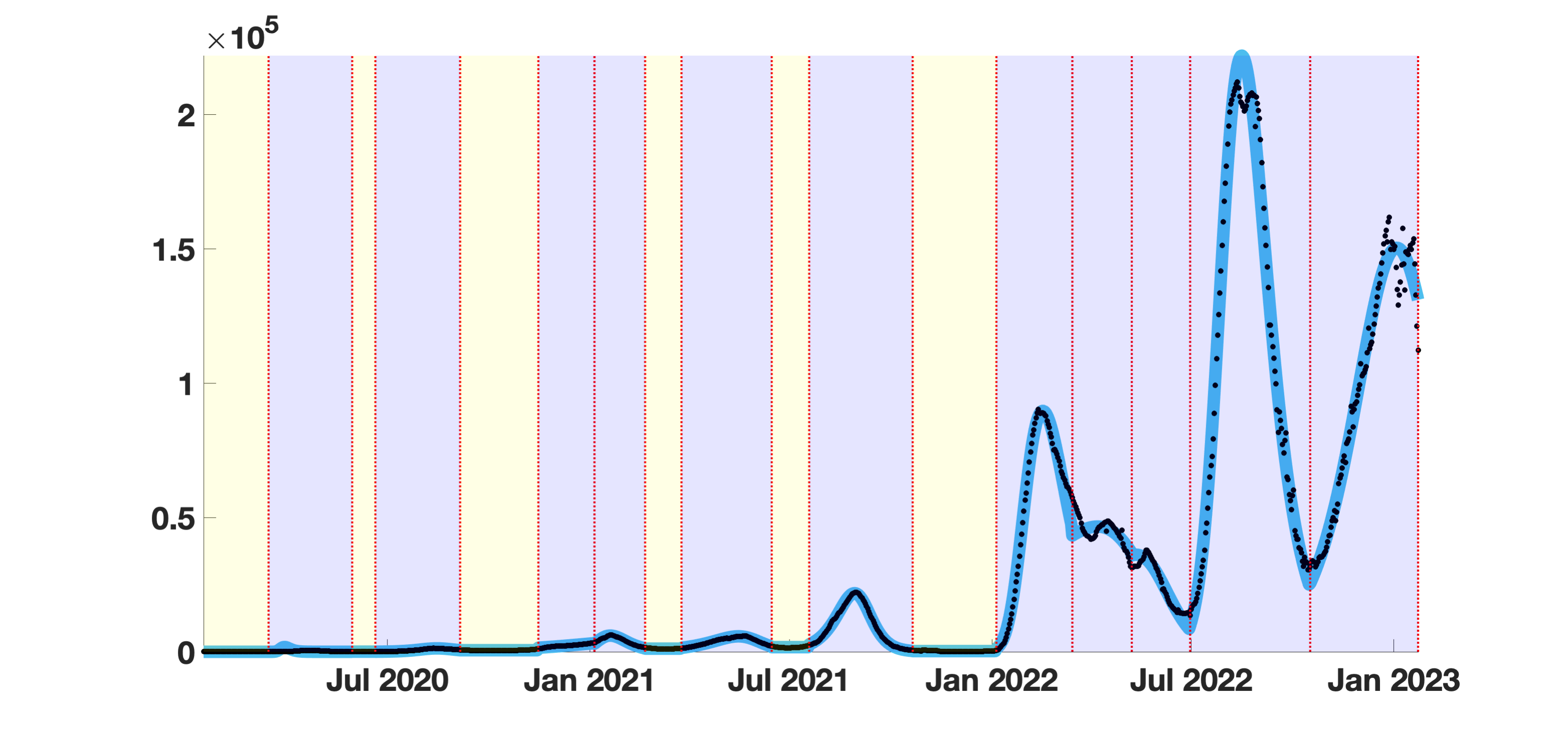}
		
	\end{center}
	\caption{\textit{In this figure we plot   in blue the first derivative of the phenomenological model  and in black the data. Data is the daily  reported  number of  new cases  with a 14-day rolling average. }}\label{Fig11}
\end{figure}

%
%
%
%
%
%
%
%
%
%
%
By making a principal component analysis (i.e. the matlab function PCA) between the standardized variables $CV_s(N_t), Skew_s(N_t),  Kurt_s(N_t)$, and $ \mathcal{E}_s(N_t)$, we obtain the percentage of the variance explained  by each principal component
$$
Explain= \left( \begin{array}{cccc}
 71.62\\
17.54 \\
6.87\\
3.97
\end{array}
\right)
$$ 
and the matrix giving the projection coefficients of the principal components 
$$
coeff =\left( \begin{array}{cccc}
 0.5234 &  -0.2243 &   0.7969  & -0.2018\\
0.5452   &-0.2406  & -0.2310 &   0.7691\\
0.5401  & -0.1760 & -0.5575 &  -0.6054\\
0.3703  &  0.9278   & 0.0270   & 0.0358

\end{array}
\right).
$$  
By using the first column of the above matrix, we deduce the first principal component 
\begin{equation} \label{10}
	0.52 CV_s(N_t)+  0.55  Skew_s(N_t)+ 0.54 Kurt_s(N_t)+ 0.37 \mathcal{E}_s(N_t)
\end{equation}
which explains $71 \%$ of the variability.

We deduce that Kurtosis, Skewness, and the coefficient of variation (in decreasing order of importance) best explain the variability. 

%

\begin{figure}[H]
	\begin{center}
		
		\includegraphics[scale=0.25]{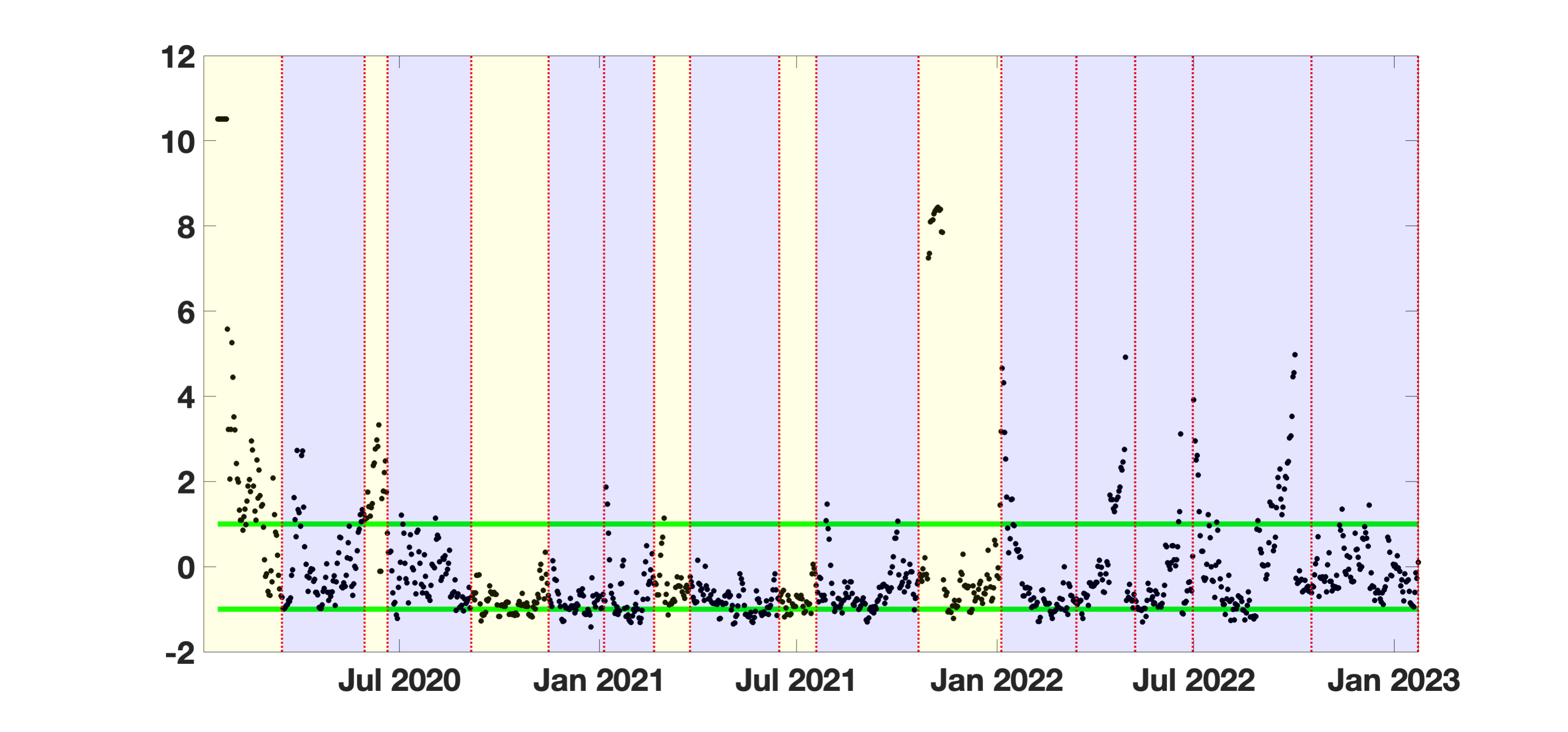}
		
	\end{center}
	\caption{\textit{In this figure we plot  the first principal component for Japan (see formula \eqref{10}). The horizontal green lines correspond to the values $\pm 1$.    }}\label{Fig12}
\end{figure}

\section{Discussion}
The forecasting of the epidemic waves of the COVID-19 outbreak is based on a change in the nature of the time series dynamics related to the number of daily new reported cases of this contagious disease. This change can concern the moments or the entropy of the empirical distribution of the stationary component at the end of the endemic phase, which disrupts when a not constant trend occurs, marking the start of an epidemic wave. 

From a careful examination of Figures \ref{Fig4}, \ref{Fig8} and \ref{Fig12}, we can conclude that there are not constant but frequent patterns for $C_1(t)$ identifiable in the three studied countries and for a majority of their endemic/epidemic transitions.

The predictive power of the first principal component $C_1(t)$ can be quantified by its performance ratio, that is, by the percentage of correct retro-predictions obtained by fixing variation thresholds to forecast the occurrence of an epidemic wave. For France, if we fix the threshold to the value $1$, $C_1(t)$ predicts correctly an epidemic outbreak a weak after a decrease from this threshold value not reached elsewhere in the endemic phase. This prediction is correct at $53\%$ only for India. For Japan, the performance of the retro-prediction is $71\%$ and $70\%$ for France.

It is clear that the level of prediction is not very high (71\% in the best case), but, in the absence of a currently reliable predictor, we can consider that it is sufficient to trigger mitigation measures at level of a population. A more systematic study of the changing shape of the empirical distribution is needed, looking at many epidemic waves in many countries. A parameter measuring the deviation from classical laws (such as those linked to the Kolmogorov-Smirnov, chi-square or Shapiro-Wilk tests in the case of the normal law) could thus be added in subsequent works.

As previously noticed in \cite{Weitz}, a classical epidemic peak with a near-symmetric growth and decline may be preceded or followed by a shoulder-like behavior corresponding to a prematurely stopped wave followed by another. Shoulder-like behavior for epidemic waves can be explained by using multiple sub-group epidemic models. This idea was first explored for SARS-CoV-1 in \cite{Magal} and reconsidered for SARS-CoV-2 by \cite{Chowell}. But the changes between endemic and epidemic are still challenging to model.  
\section{Conclusions}
	
We have studied in this article the evolution of four parameters related to the dynamics of the number of the daily reported new cases $N_t$ at day $t$ of a contagious disease, which can serve as early indicators of the appearance of epidemic waves from a previous endemic state. By applying this parameter calculation to COVID-19, we showed that a score obtained by PCA based on the linear combination of the four chosen parameters with specific coefficients for each one could forecast the variations of the empirical distribution of the daily reported number of new cases $N_t$, then could be considered often as a good predictor (for the countries and the epidemic waves in these countries) of the endemic-epidemic transition. A systematic study of contagious diseases other than COVID-19 is necessary to confirm this forecasting property's existence. Still, we can already propose this score as a realistic indicator of the next occurrence of an epidemic outbreak from a change in the dynamics of the observed daily new cases during the endemic periods.

\end{document}